\documentclass{amsart}
\def\Q{{\mathbf Q}}
\def\Z{{\mathbf Z}}
\def\C{{\mathbf C}}

\def\F{{\mathbf F}}
\def\SS{{\mathbf S}}
\def\A{{\mathbf A}}
\def\Sn{{\mathbf S}_n}
\def\An{{\mathbf A}_n}

\def\Gal{\mathrm{Gal}}
\def\Perm{\mathrm{Perm}}

\def\End{\mathrm{End}}
\def\Aut{\mathrm{Aut}}

\def\Mat{\mathrm{Mat}}
\def\cl{\mathrm{cl}}
\def\supp{\mathrm{supp}}

\def\I{\mathrm{Id}}

\def\fchar{\mathrm{char}}

\def\GL{\mathrm{GL}}

\def\dim{\mathrm{dim}}

\def\P{{\mathbf P}}

\newtheorem{thm}{Theorem}[section]
\newtheorem{lem}[thm]{Lemma}

\theoremstyle{definition}

\newtheorem{exs}[thm]{Examples}

\newtheorem{rem}[thm]{Remark}

\title{Hyperelliptic jacobians without complex multiplication}
\author[Yuri G. Zarhin]{Yuri G. Zarhin}

\begin{document}
\maketitle
\section{Introduction}
The aim of this note is to prove that
in characteristic $0$ the
jacobian $J(C)=J(C_f)$  of a hyperelliptic curve
$$C=C_f:y^2=f(x)$$
has only trivial endomorphisms
over an algebraic closure of the ground field $K$
if the Galois group $\Gal(f)$ of the polynomial
$f \in K[x]$ is ``very big".

More precisely, if $f$ is a polynomial of degree $n \ge 5$
and $\Gal(f)$ is either the symmetric group $\Sn$ or the
alternating group $\An$ then $\End(J(C))=\Z$. 
Notice that it easily
follows that  the ring of $K$-endomorphisms of $J(C)$
coincides with $\Z$ and the real problem is how to prove that
every endomorphism of $J(C)$ is defined over $K$.
 
There some results of this type in the literature.
Previously Mori \cite{Mori1}, \cite{Mori2} has constructed
explicit examples (in all characteristic)
of hyperelliptic jacobians without nontrivial endomorphisms.
In particular, he provided examples over $\Q$ with
semistable $C_f$ and big (doubly transitive) $\Gal(f)$
\cite{Mori2}. The semistability of $C_f$ implies the
semistability of $J(C_f)$ and, thanks to a theorem of Ribet 
 \cite{Ribet}, all endomorphisms of $J(C_f)$ are defined
over $\Q$. 
(Applying to $C_f/\Q$ the Shafarevich conjecture \cite{Sh} (proven by Fontaine \cite{F} and independently by Abrashkin \cite{A1}, \cite{A2}) and
using Lemma 4.4.3 and arguments on p. 42 of \cite{Serre}, 
one may prove that the Galois group $\Gal(f)$
of the polynomial $f$ involved is $\SS_{2g+1}$ where
$\deg(f)=2g+1$.)

Andr\'e (\cite{Masser}, pp. 294-295) observed that results of Katz (\cite{Katz1}, \cite{Katz2})
give rise to examples of hyperelliptic jacobians $J(C_f)$ over the field of rational function $\C(z)$ with $\End(J(C_f))=\Z$. Namely, one may take $f(x)=h(x)-z$ where $h(x) \in \C[x]$ is a
{\sl Morse function}. In particular, this explains Mori's example \cite{Mori1} 
$$y^2=x^{2g+1}-x+z$$
over $\C(z)$. 

Notice  that if $h$ is a {\sl Morse polynomial} of degree $n$ then the Galois group of $h(x)-z$ over $\C(z)$ is the symmetric group $\SS_n$ (\cite{Serre}, Th. 4.4.5, p. 41).

Masser \cite{Masser} constructed a completely different class of hyperelliptic jacobians $J(C_f)$ over $\C(z)$ with $\End(J(C_f))=\Z$.  In his examples $f$ splits into a product of linear factors over $\C(z)$ as follows.
  $$f(x)=x(x-z^{\alpha(1)})(x-z^{\beta(1)}) \dots (x-z^{\alpha(g)})(x-z^{\beta(g)})$$ 
where 
$$0 \le \alpha(1)<\beta(1)< \cdots < \alpha(g) < \beta(g)$$
 and the differences $\alpha(1)-\beta(1), \dots , \alpha(g)-\beta(g)$ are distinct. Masser's proof is purely analytic in character.

This paper was written during my stay in Glasgow. I would like
to thank the Department of Mathematics of University of Glasgow
for its hospitality.

\section{Main result}
\label{mainr}
Throughout this paper we assume that $K$ is a field of characteristic
different from $2$. We fix its algebraic closure $K_a$
and write $\Gal(K)$ for the absolute Galois group $\Aut(K_a/K)$.

\begin{thm}
\label{main}
Let $K$ be a field with $\fchar(K) \ne 2$,
 $K_a$ its algebraic closure,
$f(x) \in K[x]$ an irreducible polynomial of
degree $n \ge 5$ such that the Galois group of $f$
is either $\Sn$ or $\An$. Let $C_f$ be a hyperelliptic
curve $y^2=f(x)$. Let  $J(C_f)$ be  its jacobian, $\End(J(C_f)))$ the 
ring of $K_a$-endomorphisms of $J(C_f)$. 
Then either $\End(J(C_f))=\Z$ or 
$\fchar(K)>0$ and $J(C_f)$ is
a supersingular abelian variety.
\end{thm}

\begin{exs}
\begin{enumerate}

\item
The polynomial $x^n-x-1$ has Galois group $\Sn$ over $\Q$ (\cite{Serre}, p. 42). Hence 
 the jacobian of the curve $y^2=x^n-x-1$ has no nontrivial 
endomorphisms over $\bar{\Q}$ and therefore over $\C$ 
for all $n \ge 5$.
\item
The  Galois group of the ``truncated exponential"
 $$\exp_n(x):=1+x+\frac{x^2}{2}+\frac{x^3}{6}+ \cdots +
\frac{x^n}{n!}\in \Q[x]$$
is either $\Sn$ or $\An$ \cite{Schur}.
Hence the jacobian of the curve $y^2=\exp_n(x)$ has no nontrivial 
endomorphisms over $\bar{\Q}$ and therefore over $\C$ 
for all $n \ge 5$.

\end{enumerate}
\end{exs}

\begin{rem}
\label{odd}
Let 
$f(x) \in K[x]$ be an irreducible polynomial of even
degree $n=2m \ge 5$ such that the Galois group of $f$
is either $\Sn$ or $\An$. Then $n \ge 6$. Let
$\alpha \in K_a$ be a root of $f$ and $K_1=K(\alpha)$
be the corresponding subfield of $K_a$. We have
$$f(x)=(x-\alpha) f_1(x)$$
with $f_1(x) \in K_1[x]$. Clearly, $f_1(x)$
 is an irreducible
polynomial over $K_1$ of odd degree $2m-1=n-1 \ge 5$,
whose Galois group is either $\SS_{n-1}$ or $\A_{n-1}$
respectively.
It is also clear that the polynomials
$$h(x)=f_1(x+\alpha), h_1(x)=x^{n-1}h(1/x) \in K_1[x]$$
are irreducible of odd degree $2m-1=n-1 \ge 5$ with
the same Galois group equal $\SS_{n-1}$ or $\A_{n-1}$ 
respectively.

The standard substitution
$$x_1=1/(x-\alpha), y_1=y/(x-\alpha)^m$$
establishes a birational isomorphism between $C_f$ and
a hyperelliptic curve
$$C_{h_1}: y_1^2=h_1(x_1).$$
It follows readily that in order to prove
 Theorem \ref{main} it suffices to do the case 
of {\sl odd} $n$.
\end{rem}

We deduce Theorem \ref{main} from the following auxiliary statement.

\begin{thm}
\label{main2}
Suppose   $n=2g+1$ is an odd integer which is greater or equal than
$5$. Suppose $f(x) \in K[x]$ is a polynomial of degree $n$, whose Galois
group is either $\An$ or $\Sn$. Suppose $C$ is a hyperelliptic curve
$y^2=f(x)$ of genus $g$ over $K$, Suppose 
 $J(C)$  is the jacobian of $C$ and $J(C)_2$ is
the group of is points of order $2$,  viewed as a $2g$-dimensional
$\F_2$-vector space provided with the natural action of $\Gal(K)$.

Let $R$ be a subalgebra of $\End_{\F_2}(J(C)_2)$ which
 contains the identity operator $\I$. 
Assume that for each $u \in R, \sigma \in \Gal(K)$
 the subalgebra $R$ contains
$$^{\sigma}u: x \mapsto \sigma u \sigma^{-1}(x),
 \quad x \in J(C)_2.$$
Either $R=\F_2 \cdot \I$ or $R=\End_{\F_2}(J(C)_2)$.
\end{thm}

We prove Theorem \ref{main2} in Section \ref{main2p}.
In the next section we deduce Theorem \ref{main} from 
Theorem \ref{main2}.

\section{Proof of main result}
So, we assume that $f(x) \in K[x]$ satisfies the conditions
of Theorem \ref{main}. In light of Remark \ref{odd}, we may assume
that $n=2g+1$ is odd. Therefore $J(C)$ is the $g$-dimensional
abelian variety defined over $K$.

Since $J(C)$ is defined over $K$, one may associate with every
$u \in \End(J(C), \sigma \in \Gal(K)$ an endomorphism
$^{\sigma}u \in \End(J(C))$ such that
$$^{\sigma}u(x)=\sigma u(\sigma^{-1}x) \quad 
\forall x \in J(C)(K_a).$$ 
Let us put
$$R:=\End(J(C)) \otimes \Z/2\Z \subset \End_{\F_2}(J(C)_2).$$
Clearly, $R$ satisfies all the conditions of Theorem \ref{main2}.
This implies that either $R=\F_2 \cdot \I$ or
 $R=\End_{\F_2}(J(C)_2)$.
If $\End(J(C)) \otimes \Z/2\Z =R=\F_2 \cdot \I$
 then the free abelian group $\End(J(C))$
has rank $1$ and therefore coincides with $\Z$.
If $\End(J(C() \otimes \Z/2\Z =R=\End_{\F_2}(J(C)_2)$ 
then the free abelian group $\End(J(C))$
has rank $(2\dim(J(C)))^2=(2g)^2$ and therefore
 the semisimple $\Q$-algebra
$\End^0(J(C))=\End(J(C))\otimes\Q$ has dimension $(2g)^2$.

Now, Theorem  \ref{main} becomes an immediate corollary of the following assertion.

\begin{lem}
Let $X$ be an abelian variety of dimension $g$ over an algebraically closed field $F$. Assume that  the semisimple $\Q$-algebra
$\End^0(X)=\End(X)\otimes\Q$ has dimension $(2g)^2$.
Then $\fchar(F)>0$ and $X$ is supersingular.
\end{lem}

\begin{proof}
Let us fix a prime $\ell \ne \fchar(F)$ and consider
 the $\ell$-adic Tate module $T_{\ell}(X)$ of $X$.  Let
$V_{\ell}(X)=T_{\ell}(X)\otimes_{\Z_{\ell}}\Q_{\ell}$ be the corresponding
$\Q_{\ell}$-vector space of dimension $2g$. There is a canonical
embedding
$$\End^0(X)\otimes_{\Q}\Q_{\ell} \hookrightarrow
 \End_{\Q_{\ell}}(V_{\ell}(X))$$
and dimension arguments imply that this embedding is an isomorphism.
In particular, $\End^0(X)\otimes_{\Q}\Q_{\ell}$ is isomorphic
to the matrix algebra of size $2g$ over $\Q_{\ell}$.
Since the center of the matrix algebra over $\Q_{\ell}$ has dimension
$1$ over $\Q_{\ell}$,
the center of $\End^0(X)$ 
has dimension $1$ over $\Q$ and therefore
coincides with $\Q$.
This implies that
$\End^0(X)$ is a central simple $\Q$-algebra
of dimension $(2g)^2$. Hence, there exists
a {\sl simple} abelian variety $Y$ over $F$ and a positive integer
$r$ such that $X$ is isogenous to $Y^r$ over $F$.
 This implies that
$$g=\dim(X)=r \dim(Y),$$
$\End^0(Y)$ is a division algebra over $\Q$ and
$\End^0(X)$ is isomorphic to the matrix algebra of size $r$
over $\End^0(Y)$. In particular, 
$$\dim_{\Q}(\End^0(X))=r^2 \dim_{\Q}(\End^0(Y)).$$

Since the center of $\End^0(X)$ coincides with
$\Q$, the center of $\End^0(Y)$ also coincides with $\Q$.
It follows from Albert's classification (\cite{MumfordAV}, Sect. 21)
that either $\End^0(Y)=\Q$ or $\End^0(Y)$ is a definite quaternion
algebra over $\Q$.

If $\End^0(Y)=\Q$ then $\End^0(X)$ has dimension
$r^2 \le
(r \dim(Y))^2=g^2< (2g)^2$.
This implies that $\End^0(Y)$ is a quaternion $\Q$-algebra
and therefore
$$\dim_{\Q}(\End^0(X))=r^2 \dim_{\Q}(\End^0(Y))=4r^2=(2r)^2.$$
On the other hand, $\dim_{\Q}(\End^0(X))=(2g)^2$. This
implies that $2r=2g$, i.e., $r=g=\dim(X)$ and $Y$ is an
elliptic curve. Since $\End^0(Y)$ is the quaternion algebra,
$Y$ is a supersingular elliptic curve and $\fchar(F) > 0$.
Since $X$ is isogenous to $Y^r$, it is a supersingular abelian
variety.
\end{proof}

\section{Points of order $2$ on hyperelliptic jacobians}
\label{main2p}

Let $C$ be a hyperelliptic curve over $K$ defined by an
equation $y^2=f(x)$ where $f(x) \in K[x]$ is a polynomial
of odd degree $n$ without multiple roots. The rational function 
$x \in K(C)$ defines a canonical double cover $\pi:C \to \P^1$. 
Let $B'\subset C(K_a)$ be the set of ramification points 
of $\pi$ (Weierstra{\ss } points). Clearly, the restriction of
$\pi$ to $B'$ is an injective map 
$\pi:B' \hookrightarrow \P^1(K_a)$, whose image is the disjoint
union of $\infty$ and the set $R_f$ of roots of $f$. By abuse
of notation, we also denote by $\infty$ the ramification point
lying above $\infty$. Clearly, $\infty \in C(K)$. 
We denote by $B$ the complement of $\infty$ in $B'$.
Clearly,
$$B=\{(\alpha,0)\mid f(\alpha)=0\} \in C(K_a)\}$$
and $\pi$ defines a bijection between $B$ and $R_f$ which
commutes with the action of  $\Gal(K)$.

We write $Q_B$ for the $\F_2$-vector space of
  subsets in $B$ of {\sl even} cardinality with symmetric difference
as a sum.  There is a natural structure of $\Gal(K)$-module on $Q_B$.

Here is an explicit description of the group $J(C)_2$ of points
of order $2$ on the jacobian $J(C)$. 
Namely, let $T \subset B'$ be a subset of even cardinality.
Then (\cite{Mumford}, Ch. IIIa, Sect. 2, Lemma 2.4); see also
\cite{Mori2})
the divisor
$e_T=\sum_{P \in T}(P) -\#(T)(\infty)$ on $C$ has degree $0$
and $2 e_T$ is principal. If $T_1,T_2$ are two subsets of even
cardinality in $B'$ then the divisors $e_{T_1}$ and $e_{T_2}$
are linearly equivalent if and only if either $T_1=T_2$ or
$T_2=B'\setminus T_1$. Also, if $T=T_1\triangle T_2$ then
the divisor $e_T$ is linearly equivalent to $e_{T_1}+e_{T_2}$.

Counting arguments imply easily that each point of $J(C)_2$
is the class of $e_T$ for some $T$. We know that such a choice
is not unique. However, if we demand that $T$ does not contain
$\infty$ then such a choice always exists and unique. This
observation leads to a canonical group isomorphism
$$Q_B \cong J(C)_2, \quad T \mapsto \cl(e_T).$$
Here $\cl$ stands for the linear equivalence class of a divisor.
Clearly, this isomorphism commutes with natural actions of
$\Gal(K)$. In other words, the $\Gal(K)$-modules $Q_B$
and $J(C)_2$ are canonically isomorphic.

One may describe explicitly the Galois action on $Q_B$.
In order to do  that let us consider the splitting field
$L \subset K_a$ of $f$  and let 
$G=\Gal(L/K)$ be its Galois group. Clearly, $G$ may be viewed
as a group of permutations of $R_f$ and therefore (via $\pi$)
as a subgroup in the group $\Perm(B)$ of permutations of $B$.  This
induces  obvious embeddings
$$G \subset \Perm(B) \subset \Aut(Q_B)$$
and $\Gal(K)$ acts on $Q_B$ via the composition of
the canonical {\sl surjection} $\Gal(K) \to \Gal(L/K)=G$
and the embedding 
$$G \subset \Perm(B) \subset \Aut(Q_B).$$

Now, one may easily check that Theorem \ref{main2} follows
readily from the following purely group-theoretic statement.

\begin{thm}
\label{group}
Let $B$ be a finite set of odd cardinality $n \ge 5$, 
$Q_B$ the $\F_2$-vector space of
its  subsets of {\sl even} cardinality with symmetric difference
as a sum. Let $S=\Perm(B)$ be the group of permutation of $B$ viewed as
a subgroup of $\Aut(Q_B)$.
Let $G$ be a subgroup of $S$ which is isomorphic either to $\Sn$
or to $\An$.

Let $R$ be a subalgebra of $\End_{\F_2}(Q_B)$
 which contains the identity operator $\I$. Assume that
$$u R u^{-1} \subset R \quad
 \forall u \in G \subset S \subset \Aut(Q_B).$$
Either $R=\F_2 \cdot \I$ or $R=\End_{\F_2}(Q_B)$.
\end{thm}

We prove Theorem \ref{group} in the next Section.

\section{Representation theory}
We keep all the notations and assumptions of Theorem \ref{group}..
Clearly, $S \cong \Sn$. We write $A$ for the only subgroup
in $G$ of index $2$. Clearly, $A$ is normal and isomorphic
to the alternating group $\An$. 
It is well-known that the group $A$ is simple
 of order $n!/2$.
The cardinality arguments
imply  easily that either $G=S$ or $G=A$.

We have
$$A \subset S \subset \Aut(Q_B), \quad
\dim_{\F_2}(Q_B)=n-1.$$
 
Let us consider the $n$-dimensional $\F_2$-vector space $\F_2^B$ of 
all subsets of $B$ with symmetric difference
as a sum.  The space $\F_2^B$ is  provided with a natural action of  $S=\Perm(B)$. 
It is well-known that one may view $\F_2^B$ as the $\F_2$-vector space of all maps $f: B \to \F_2$ provided with a natural action of $S$. Namely, a subset $T$ corresponds to its characteristic function $\chi_T:B \to \{0,1\}=\F_2$ and  a function $\varphi: B \to \F_2$ corresponds to its support
$\supp(\varphi)=\{x  \in B\mid \varphi(x)=1\}$.

Cl;early, $Q_B \subset \F_2^B$ is a $S$-stable hyperplane in $ \F_2^B$
and two actions of $S$ on $Q_B$ do coincide. 
 Since $n=\#(B)$ is odd, the set $B$ does not belong to $Q_B$. Clearly, $B \in \F_2^B$ is $S$-invariant.
 Let $L$ be the one-dimensional subspace of $\F_2^B$ generated by $B$. Clearly, $S$ acts trivially on $L$ 
  and  there  is a $S$-invariant splitting
$$\F_2^B=Q_B \oplus L$$
which is also $A$-invariant. It is also clear that $\End_A(L)=\End_{\F_2}(L)=\F_2$.
This implies that $\dim_{\F_2}(\End_A(L))=1$. It is also clear that $\dim_{\F_2}(\End_A(Q_B)) \ge 1$.

\begin{lem}
\label{absirr}
$\End_A(Q_B)=\F_2$.
\end{lem}

\begin{proof}
It suffices to check that $ \dim_{\F_2}(\End_A(Q_B)) \le 1$.

The $A$-invariant splitting $\F_2^B=Q_B \oplus L$ implies that
$$\dim_{\F_2}(\End_A(\F_2^B)) \ge \dim_{\F_2}(\End_A(Q_B))+ \dim_{\F_2}(\End_A(L))=$$ 
$$\dim_{\F_2}(\End_A(Q_B))+ 1.$$
This implies that $\dim_{\F_2}(\End_A(\F_2^B)) \ge \dim_{\F_2}(\End_A(Q_B))+ 1$. 
Since $A$ acts doubly transitively on $B$, we have
$ \dim_{\F_2}(\End_A(\F_2 ^B))=2$ (\cite{Passman}, Lemma 7.1 on p. 52). This implies that $1 \ge \dim_{\F_2}(\End_A(Q_B))$.
\end{proof}

\begin{lem}
\label{an}
The $A$-module $Q_B$ is absolutely simple.
\end{lem}

\begin{proof}
Let us prove first that $Q_B$ is irreducible.
Let $U$ be a non-zero $A$-stable
subspace in $Q_B$. Let $T \in U$ be a non-empty subset
 of $B$ with smallest  possible cardinality.
Since $n$ is odd, $T \ne B$.
 If $T$ consists of $2$ elements
then we are done, because $A=\An$ acts
 doubly transitively on $B$ and each
subset in $B$ of even cardinality could be presented as a
symmetric difference (disjoint union) of
$2$-element sets. So, assume that $T$ consists
 of at least $4$ elements.
Pick elements $t \in T$ and $b \in B \setminus T$. 
Then there is an {\sl even}
permutation $s \in A$ such that $s(T)=T\setminus \{t\} \bigcup \{b\}$.
Clearly, the symmetric difference $T \triangle s(T) \in U$ has two elements
less than $T$ which contradicts the choice of $T$. This proves the
irreducibility of $Q_B$.

It follows from Lemma \ref{absirr} that $Q_B$ is {\sl absolute} irreducible.
\end{proof}

Since $G$ always contains $A$, Theorem \ref{group} is an immediate
corollary of the following statement.

\begin{thm}
\label{An}
Let $R$ be a subalgebra of $\End_{\F_2}(Q_B)$ 
which contains the identity operator
$\I$. Assume that
$$u R u^{-1} \subset R \quad \forall u \in A \subset \Aut(Q_B).$$
Either $R=\F_2 \cdot \I$ or $R=\End(Q_B)$.
\end{thm}

\begin{proof}[Proof of Theorem \ref{An}]
Clearly, $Q_B$ is a faithful $R$-module and
$$u R u^{-1} = R \quad \forall u \in A \subset \Aut(Q_B).$$

{\bf Step 1}. $Q_B$ is a {\sl semisimple} $R$-module. Indeed, 
let $U \subset Q_B$
be a simple $R$-submodule. Then
$U'=\sum_{s\in A} sU$ is a non-zero $A$-stable subspace in $Q_B$
and therefore must coincide with $Q_B$. 
On the other hand, each $sU$ is
also a  $R$-submodule in $Q_B$, because $s^{-1} R s=R$. In addition,
if $W \subset sU$ is an $R$-submodule then $s^{-1}W$ is an $R$-submodule in $U$,
because 
$$R  s^{-1}W=s^{-1} s R  s^{-1}W=s^{-1} RW=s^{-1}W.$$
 Since $U$ is simple, $s^{-1}W=\{0\}$ or $U$. This implies that $sU$ is also simple. Hence
 $Q_B=U'$ is a sum of simple $R$-modules and therefore  is a  semisimple $R$-module.

{\bf Step 2}. The $R$-module $Q_B$ is {\sl isotypic}. Indeed, let us split the semisimple
$R$-module $Q_B$ into the direct sum
$$Q_B =V_1 \oplus \cdots \oplus V_r$$
of its isotypic components.
 Dimension arguments imply that $r \le \dim(Q_B) = n-1$.
 It follows easily from the arguments of the previous step that for each isotypic component $V_i$ its image 
  $s V_i$ is an isotypic $R$-submodule for each $s \in A$ and therefore is contained in some $V_j$.
  Similarly, $s^{-1}V_j$ is an isotypic submodule obviously containing $V_i$. Since $V_i$ is the isotypic component,
  $s^{-1}V_j=V_i$ and therefore $sV_i=V_j$.
  This means that $s$ permutes the $V_i$; since $Q_B$ is $A$-simple, $A$ permutes them transitively.
 This gives rise to the homomorphism $A \to \SS_r$ which must be trivial,  
since $A$ is the
simple group, whose order ($=n!/2$)
is greater than $(n-1)! \ge r!=$ order of $\SS_r$. 
This means that 
$s V_i=V_i$ for all $s \in A$ and $Q_B=V_i$ is isotypic.

{\bf Step 3}. Since $Q_B$ is isotypic, there exist a simple
 $R$-module $W$
and a positive integer $d$ such that $Q_B \cong W^d$.
 Clearly,
$$d \cdot \dim(W)=\dim(Q_B)=n-1.$$
It is also clear that
$\End_R(Q_B)$ is isomorphic to the matrix algebra
 $\Mat_d(\End_R(W))$ of size $d$ over $\End_R(W)$.

Let us put
$$k=\End_R(W).$$
Since $W$ is simple, $k$ is a finite division algebra of
characteristic $2$. Therefore $k$ is a finite field
 of characteristic $2$ and $[k:\F_2]$ divides $n-1$.
We have $\End_R(Q_B) \cong \Mat_d(k)$.
CLearly, $\End_R(Q_B) \subset \End_{\F_2}(Q_B)$ is stable under the adjoint
action of $A$. This induces a homomorphism
$$\alpha: A \to \Aut(\End_R(Q_B))=\Aut(\Mat_d(k)).$$
Since $k$ is the center of $\Mat_d(k)$,
 it is stable under the action of $A$,
i.e., we get a homorphism
$A \to \Aut(k)$, which must be trivial,
 since $A$ is the simple group and $\Aut(k)=\Gal(k/\F_2)$ is abelian.
 This implies that
the center $k$ of $\End_R(Q_B)$ commutes with $A$. 
Since $\End_A(Q_B)=\F_2$, we have $k=\F_2$.   
This implies that
$\End_R(Q_B) \cong \Mat_d(\F_2)$ and 
$$\alpha: A \to \Aut(\Mat_d(\F_2))=\GL(d,\F_2)/\F_2^*=\GL(d,\F_2)$$
 is trivial if and only if $\End_R(Q_B) \subset \End_A(Q_B)=\F_2 \cdot \I$.
 Since $\End_R(Q_B) \cong \Mat_d(\F_2)$, $\alpha$ is trivial if and only if
$d=1$, i.e. $V$ is an absolutely simple $R$-module.
 It follows from Jacobson density
theorem that
$R \cong \Mat_m(\F_2)$ with $dm=n-1$.
This implies that $\alpha$ is trivial if and only if
$R \cong \Mat_{n-1}(\F_2)$, i.e., $R=\End(Q_B)$.

The adjoint action of $A$ on $R$ gives rise
 to a homomorphism
$$\beta: A \to \Aut(\Mat_m(\F_2))=\GL(m,\F_2).$$
Clearly, $\beta$ is trivial if and only if 
$R$ commutes with $A$,
i.e. $R=\F_2 \cdot \I$.

{\bf Step 4}. It follows from the previous step that we are done if either
$\alpha$ or $\beta$ is trivial. Since $md = n-1$ , 
the desired triviality
follows from the following 

{\bf Observation}. Let $c$ be a positive integer. 
If $c^2 \le n-1$ then
every homomorphism from $A \cong A_n$ to $\GL(c, \F_2)$
is trivial. Indeed, the cardinality of $\GL(c,\F_2)$ is
 strictly less than $2^{c^2}$ and
$2^{c^2} \le 2^{n-1}<n!/2=$ order of the {\sl simple} group $A$ ($n \ge 5$).

\end{proof}


\bigskip
 
\noindent {\small {Department of Mathematics, Pennsylvania State University,}}
 
\noindent {\small {University Park, PA 16802, USA} }
 
\vskip .4cm

\noindent {\small {Institute for Mathematical Problems in Biology,}}

\noindent {\small {Russian Academy of Sciences, Push\-chino, Moscow Region, 142292, RUSSIA}}
 
 \vskip .4cm
 
\noindent {\small  {Department of Mathematics, University of Glasgow,}}

\noindent {\small  {15 University Gardens, Glasgow G12 8QW, Scotland, UK}}
 
 \vskip  .4cm

\noindent {\small {\em E-mail address}: zarhin@math.psu.edu}

\end{document}